\input amstex
\documentstyle{amsppt}
\magnification=1200
\overfullrule=0pt
\NoRunningHeads
%
\def\C{\Bbb C}

\def\R{\Bbb R}

\def\Z{\Bbb Z}

%
\topmatter
\title
Continuous cohomology of the group of volume-preserving and symplectic diffeomorphisms, measurable transfer and higher asymptotic cycles
\endtitle
%
\author
Alexander Reznikov
\endauthor
%
%
%
%
\address
Institute of Mathematics, Hebrew University, Giv'at Ram 91904, Jerusalem, Israel, simplex$\@$math.huji.ac.il
\endaddress
\curraddr
Max-Planck-Institut f\"ur Mathematik, Gottfried-Claren-Str. 26, 53225 Bonn, Germany, reznikov$\@$mpim-bonn.mpg.de
\endcurraddr
%
%
\date
July, 1996
\enddate
%
%
%
%
\thanks
Partially supported by a BSF grant
\endthanks
%
%
%
%
%
%
%
%
%
%

\endtopmatter
\document

Topology of a manifold is reflected in its diffeomorphism group. It is challenging therefore to understand the diffeomorphism group
$Diff(M)$
both as a topological and discrete group. Twenty years ago, some work has been done, in connection with characteristic classes of foliations, in constructing continuous cohomology classes for
$Diff(M)$.
For $M$ closed oriented $n$-dimensional manifold, a class in
$H^{n+1}_{cont}(Diff(M),\R)$
has been explicitly written down by Bott [Bo] [Br]. This class is defined as follows. The group
$Diff(M)$
acts in the multiplicative group
$C^{\infty}_+(M)$
of positive smooth functions, and on its torsor
$A_n(M)$
of volume forms. Hence one gets a cocycle in
$H^1_{cont}(Diff(M),C^{\infty}_+(M))$,
defined by
$\lambda(f) = \frac{f^*(v)}v = Jac_v(f)$,
where
$\nu \in A_n(M)$
and
$f \in Diff(M)$.
The Bott class is

$$\int_M \log \lambda \cup \undersetbrace n \to {d \log \lambda \cup \ldots \cup d \, \log \lambda}$$
The nontriviality of Bott class had been shown for
$M = S^1$
[Br], and recently for
$S^n$
[BCG],
$\C P^n$
[Go] by restricting to finite-dimensional Lie groups in
$Diff(M)$.
In fact, the restriction of the Bott class on
$SO(n,1) \subset Diff(S^n)$
gives the hyperbolic volume class, whereas the restriction on
$PSL(n+1, \C) \subset Diff (\C P^n)$
gives the Borel class.

By its construction, the Bott class vanishes on the group
$Diff_{\nu}(M)$
of volume-preserving diffeomorphisms. Moreover, since it is defined by an invariant closed
$(n+1)$
-form in the space
$A_n(M)$
where
$Diff(M)$
acts, and by a theorem of Brooks [Br] there are no more invariant forms there, one gets just one class in dimension
$(n+1)$
for a fixed manifold $M$. This contrasts sharply the usual intuition coming from the study of finite-dimensional semisimple group, where there is a range of continuous cohomology classes.

In this paper we construct, for a closed manifold
$M^n$
with a volume form
$\nu$, 
a series of continuous cohomology classes in
$H^{\kappa}_{cont}(Diff_{\nu}(M),\R)$
for all
$\kappa = 5,9, \ldots$. 
The classes will be shown nontrivial already for a torus
$T^n$.
We also will construct, for a symplectic manifold
$(M,w)$,
a series of classes in
$H^{2\kappa}(Sympl (M),\R)$
for
$\kappa = 1,3, \ldots$. 
Again, these are nontrivial for a torus
$T^n$
with standard symplectic structure.

Working harder, we will show that for the smooth moduli space of stable vector bundles over a Riemann surface
$\Cal M$
with its K\"ahler structure, our class in
$H^2(Symp(\Cal M_g),\R)$
is nontrivial and restricts to a generator of
$H^2(Map_g,\R)$,
where
$Map_g$
is the mapping class group:

\proclaim{Theorem (3.6)}
$H^2(Sympl(\Cal M_g),\R)$
is nontrivial. Moreover, the homomorphism
$Map_g \rightarrow Sympl(\Cal M_g,\R)$
induces a nontrivial map in the second real cohomology.
\endproclaim

In both cases, our classes arise from action on a ``principal homogeneous space'' $X$ which in the case of
$Diff_{\nu}(M)$
will be the space of Riemannian metrics with volume form
$\nu$, 
and in the case of
$Sympl(M)$
will be the twistor variety, introduced in [Re1]. In that paper we have studied the symplectic reduction of $X$ with respect to the Hamiltonian action of subgroups of
$Sympl(M)$
with a primal interest in integrable systems arising on Teichm\"uller space and universal Jacobian. A lenghty computation from [Re1] related to the existence of the moment map will be used here to prove a vanishing result in 5.4.

There is quite another way to look at our classes, from the stand point of the transfer map. The subgroup
$Diff^0_{\nu}(M)$
of
$Diff_{\nu}(M)$
which fixes a point
$p \in M$,
has the tangential representation to
$SL_n(\R)$
and one can pull the Borel classes back on
$Diff^0_{\nu}(M)$.
The transfer map [Gu] will send these classes to
$H^*_{cont}(Diff_{\nu}(M))$.
We will not however prove a rigorous comparison theorem relating these two types of construction in the present paper. However we do use the transfer map to define a new source of classes in
$H^*(Diff_{\nu}(M))$
coming from the fundamental group of $M$. Namely, a map

$$S : H^{\kappa}(\pi_1(M),\R) \rightarrow H^{\kappa}(Diff^{\sim}_{\nu}(M),\R)$$
will be constructed where
$Diff^{\sim}_{\nu}(M)$
is the connected component of
$Diff_{\nu}(M)$.
For
$\kappa =1$,
the dual of this map, a character

$$S^{\vee} : Diff^{\sim}_{\nu} (M) \rightarrow H_1(M,\R)$$
has been known for forty years [Sch] and called the asymptotic cycle map. One can view our map $S$ as ``higher'' asymptotic cycle map.

For $M$ a closed surface with an area form, the groups
$Diff_{\nu}(M)$
and
$Sympl(M)$
coincide. The two previously described constructions produce a class in
$H^2_{cont}(Diff_{\nu}(M))$
which we will show to lie in bounded cohomology group
$H^2_b(Diff(M),\R)$.
For
$f,g \in Diff_{\nu}(M)$
we give an explicit formula for a cocycle
$\ell (f,g)$
representing this class. For any lamination on $M$ [Th] one can exhibit quite a different formula, using the expression for Euler class from [BG].

The following application of dynamical nature will be proven. Let
$F_2$
be a free group in two generators, and let, for some words
$h_i,k_i$
in
$F_2$,
a sum
$\sum^{\infty}_{i=1} a_i(h_i,k_i),$ 
$\Sigma |a_i| < \infty$
be a cycle for
$\ell^1$
-homology of
$F_2$.
This homology has dimension
$2^{\aleph_0}$,
as shown in [M]. Let $M$ be a closed surface with an area form $\nu$. Given
$f,g \in Diff_{\nu}(M)$
one has a homomorphism
$F_2 \rightarrow Diff_{\nu}(M)$,
so the words
$h_i,k_i$
may be viewed as diffeomorphisms in
$Diff_{\nu}(M)$.

\proclaim{Theorem (4.2)}
Suppose
$\sum^{\infty}_{i=1} a_i \, \ell (h_i,k_i) \ne 0$.
Then the group generated by
$f,g$
in
$Diff_{\nu}(M)$
is not amenable.
\endproclaim

The significance of Theorem 4.2 stems from the fact that the condition
$\Sigma \, a_i \ell(h_i,k_i) \ne 0$
is
$C^1$
-open on
$f,g$.
Therefore one gets a \underbar{domain} in
$Diff_{\nu}(M) \times Diff_{\nu}(M)$,
such that any pair
$(f,g)$
in it generate a ``big'' group in
$Diff_{\nu}(M)$.
One can see this result as a step towards ``Tits alternative'' for the infinite-dimensional Lie group
$Diff_{\nu}(M)$.

We will show in the next paper that this theorem holds for $M$ symplectic of higher dimension. For that purpose we ill use Lagrangian measurable foliations and Lyon-Vergne Maslov class to show that our class in
$H^2(Sympl(\Cal M_g,\R)$
is bounded. See also the end of [BG].

In [Re2] we defined the ``symplectic Chern-Simons'' classes
$K^{alg}_{2i-1}(Sympl (M)) = \pi_{2i-1} ((B \, Symp)^{\delta} (M))^+) \rightarrow \R/A$,
where $A$ is the group of periods of the Cartan form in
$\Omega^{2i-1}_{cl}(\text{Sympl}^{top}(M))$,
introduced in [Re2], on the Hurewitz image of
$\pi_{2i-1}(\text{Sympl}^{top}(M))$
in
$H_{2i-1}(\text{Sympl}^{top}(M),\R)$.
The real classes introduced in the present paper seem to be in the same relation to the symplectic Chern-Simons classes as Borel classes in
$H^*_{cont}(SL_n(K),\R)$
are to proper Chern-Simons classes
$(K = \R, \C)$.
The ``symplectic Chern-Simons classes'' of [Re2] have remarkable rigidity property: for a continuous family of representations of a f.g. group
$\Gamma$
into
$\text{Sympl}(M)$,
the pull-back of these classes are constant in
$H^*(\Gamma)$.
This contrasts strikingly the famous non-rigidity of the Bott class, proved by Thurston. In fact, Thurston exhibited a family of homomorphism
$\pi_1(S) \rightarrow Diff(S^1)$,
where $S$ is a closed surface of genus two, with varying Godbillon-Vey class (which coincides with the Bott class for
$Diff(S^1))$.

We do not know if the real classes constracted in the present paper in
$H^*(Diff_{\nu}(M))$
and
$H^*(\text{Sympl}(M))$
are rigid. However, we introduce a new ``Chern-Simons'' class in
$H^3(Diff_{\nu}(S^3),\R/\Z)$
which \underbar{is rigid} and restricts to usual Chern-Simons class on
$H^3(SO(4), \R/\Z)$.
This uses the invariant scalar product on Lie
$(Diff_{\nu}(S^3))$
in much the same way we used invariant polynomials on Lie 
$(\text{Sympl}(M))$
in [Re2].
\proclaim{6.6 Theorem (Chern-Simons class in $Diff_{\nu} (S^3)$ )}
There exists a rigid class in
$H^3(Diff_{\nu}(S^3), \R/\Z)$
whose restriction on
$SO(4) \approx S^3 \times S^3/\Z_2$
coincides with the sum of standard Chern-Simons classes. Moreover, for
$M = S^3/\Gamma$
there exists a class in
$H^3(Diff_{\nu}(M),\R/\Z)$
whose restriction on
$S^3$
is 
$|\Gamma |$
times the standard Chern-Simons class.
\endproclaim
\head{1. Forms on the space of metrics}\endhead

We work with the manifold $M$ with the fixed volume form
$\nu$.
Define the space
$\Cal P$
as the Frechet manifold of
$C^{\infty}$
-Riemannian metrics on $M$, whose volume form is
$\nu$.
Obviously,
$Diff_{\nu}(M)$
acts on
$\Cal P$.
We can look at
$\Cal P$
as a space of sections of a fibration
$\R \rightarrow M$
with a fiber
$SL_N(\R)/SO(N)$,
where
$N = \dim M$.
Clearly,
$\Cal M$
is contractible. For any 
$n = 5,9, \ldots$ 
fix the Borel form: a
$SL_N(\R)$
-invariant closed $n$-form on
$SL_N(\R)/SO(N)$,
normalized as in [Bo]. For a vector space $V$ of dimension $N$ with a volume form
$\nu$
this gives a canonical choice of a closed form on the space 
$\Cal P^V$
of Euclidean metrics on $V$ with determinant
$\nu$.
Call this form
$\psi^V_n$.
Now, we define a form on
$\Cal P$
by
$\psi_n =\int_M \psi_n^{T_xM} d\nu (x)$.
That means the following: let
$g \in \Cal P$
a Riemannian metric on $M$. Let
$h_1, \ldots, h_n \in T_g \Cal P$
be symmetric bilinear smooth 2-forms. Define
$\psi_n(h_1,\ldots, h_n) = \int_M \psi_n^{T_x(M)} (h_1(x), \ldots, h_n(x)) d\nu$.

\proclaim{Lemma (1.1)}
The form
$\psi \in \Omega^n(\Cal P)$
is closed and
$Diff_{\nu}(M)$
-invariant.
\endproclaim

\demo{Proof}
The invariance is obvious from definition. To prove the closedness, observe first that a form
$\psi_n(x_1,\ldots, x_m) (h_1, \ldots, h_n) = \sum^m_{j=1} \lambda_j \psi_n^{T_{x_j}(M)} (h_1(x_j), \ldots, h_n(x_j))$
is closed as a pull-back of a closed form under the map
$\Cal P \mapsto \prod^m_{j=1} \Cal P^{T_{x_j}(M)}$.
Now one approximates
$\psi$
by
$\psi_n(x_1,\ldots, x_m)$
to show that
$\psi$
is closed.
\enddemo

\demo{1.2 The definition of the classes}
We will now apply a general theory of regulators, as presented in [Re1], section 3. For a Frechet-Lie group
$\frak G$,
acting smoothly on a contractible smooth manifold $Y$, preserving a closed form
$\psi_n$,
this theory prescribes a class in
$H^n(\frak G^{\delta}, \R)$,
called
$r(\psi_n)$
in [Re1].
\enddemo

\demo{Definition (1.2)}
Consider the action of
$Diff_{\nu}(M)$
on the contractible manifold
$\Cal P$
with the invariant form
$\psi_n$
as above. A class
$\gamma_n \in H^n(Diff^{\delta}_{\nu}(M), \R)$
is defined as
$r(\psi_n)$.
\enddemo

\proclaim{Theorem (1.3)}
The class
$\gamma_n$
lies in the image of the natural map

$$H^n_{cont}(Diff_{\nu}(M),\R) \rightarrow H^n(Diff^{\delta}_{\nu} (M), \R).$$
\endproclaim

The proof follows from Proposition 1.3 below.

\demo{1.3 Simplices in $\Cal P$ and a Dupont-type construction}
Fix two metrics
$g_1,g_2$
in
$\Cal P$.
We can join them by a segment in two different ways. First, there is a straight line segment
$I_{g_1,g_2}(t) : t \mapsto t \cdot g_1 + (1-t)g_2$.
Second, there is a geodesic segment
$J_{g_1,g_2}(t): t \mapsto (x \mapsto c(t,g_1(x), g_2(x)))$.
Here
$t \in [0,1], x \in M, g_1(x), g_2(x) \in \Cal P^{T_x}(M)$
and
$c(t,g_1(x), g_2(x))$
is a geodesic segment in the homogeneous metric of symmetric space on
$\Cal P^{T_x(M)} \approx SL_N(\R)/SO(N)$.
Now, having $n$ metrics
$g_1, \ldots, g_n$
in
$\Cal P$
we define two singular simplices
$I_{g_1 \ldots g_n} : \sigma \rightarrow \Cal P$
and
$J_{g_1 \ldots g_n} : \sigma \rightarrow \Cal P$
by induction as a joint of
$g_1$
and
$I_{g_2,\ldots, g_n}$,
(resp.
$g_1$
and
$J_{g_2 \ldots g_n})$
using straight line segments (resp. geodesic segments, comp [Th2]).

Now fix a reference metric $g$ in 
$\Cal P$.
Define

$$\gamma^I_n(f_1,\ldots, f_n) = \int_{I(g,{f_1}_* g \ldots,{f_1 f_2 ...f_n}_* g )} \psi_n$$
and

$$\gamma^J_n(g_1, \ldots, g_n) = \int_{J(g,{f_1}_* g \ldots,{f_1 f_2 ...f_n}_* g )} \psi_n$$
\enddemo

\proclaim{Proposition (1.3)}
Both
$\gamma^I_n$
and
$\gamma^J_n$
are continuous cocycles, representing
$\gamma_n$.
\endproclaim

\demo{Proof}
The proof mimics the finite-dimensional case, cf. [Du], and is therefore omitted.
\enddemo

\head{2. Non-triviality}\endhead

We will prove that the class
$\gamma_n$
in discrete group cohomology, and consequently classes of
$\gamma^I_n$
and
$\gamma^J_n$
in continuous cohomology are non-trivial in general. For that purpose, consider a torus
$T^N = \R^N/\Z^N$
with a standard volume form
$dx_1 \ldots dx_N$.
We have an inclusion

$$SL(N,\Z) \hookrightarrow Diff_{\nu} (T^N)$$

\proclaim{Proposition (2.1)}
The class
$\gamma_n$
restricts to the Borel class in
$H^n(SL(N,\Z),\R)$
and is therefore nontrivial for $N$ big enough.
\endproclaim

\demo{Proof}
Let
$\Cal P_0$
be the space of left-invariant metrics on
$T^N$
with the determinant
$\nu$;
as a manifold,
$\Cal P_0 \approx SL_N(\R)/SO(N)$.
The embedding
$\Cal P_0 \hookrightarrow \Cal P$
is
$SL_N(\Z)$
-invariant, and the pull-back of the form
$\psi_n$
on
$\Cal P_0$
is the Borel form on
$\Cal P_0$.
Now by [Re1], section 3, 
$r(\psi_n)$
coincides with the Borel class.
\enddemo

\head{3. Cohomology of symplectic diffeomorphisms}\endhead

We will now adapt the theory for the group
$\text{Sympl}(M)$
of symplectic diffeomorphisms of a compact symplectic manifold $M$. For this purpose, we will introduce a new (
$\infty$
-dimensional) contractible manifold $Z$, on which
$Sympl (M)$
acts, preserving some differential forms of even degree.

\demo{3.1 Principal transformation space}
Let
$\frak F$
be the fibration over 
$M^{2n}$,
whose fiber over
$x\in M$
consists of complex structures in
$T_xM$,
say $J$,
such that
$\omega_x$
is
$J$-invariant and the symmetric form
$\omega (J \, \cdot, \cdot)$
is positive definite. Alternatively, 
$\frak F$
is a
$Sp(2n,\R)/U(n)$
fiber bundle over $M$, associated to the 
$Sp(2n,\R)$
-frame bundle. The principal transformation space $Z$ is defined as a space of
$C^{\infty}$
-sections of
$\frak F$.
So a point in $Z$ is just an almost-complex structure on $M$, tamed by $\omega$, in the sense of Gromov [Gr]. Since the Siegel upper half-plane
$Sp(2n,\R)/U(n)$
is contractible, the space $Z$ is contractible, too.
\enddemo

\demo{3.2 Forms on $Z$}
Fix an
$Sp(2n,\R)$
-invariant form on
$Sp(2n,\R)/U(N)$.
This induces a form
$\varphi^{T_xM}$
on
$\Cal F_x$
for each
$x \in M$
and a form

$$\varphi = \int_M \varphi^{T_xM} \cdot \omega^n$$
as in 1.1. Obviously, this form
$\varphi$
is
$\text{Sympl}(M)$
-invariant. Recall that the ring of
$Sp(2n,\R)$
-invariant forms on
$Sp(2n,\R)/U(n)$
is generated by forms in dimensions 
$2,6, \ldots$ 
[Bo].

Correspondingly, we have
$\text{Sympl}(M)$
-invariant closed forms, in same dimensions.
\enddemo

We single out the symplectic (K\"ahler) form on
$Sp(2n,\R)/U(n)$,
which may be described as follows. For
$J \in Sp(2n,\R)/U(n)$,
the tangent space
$T_J Sp(2n,\R)/U(n)$
consists of operators
$A : \R^{2n} \rightarrow \R^{2n}$
satisfying
$AJ = -JA$
and
$\langle Ax,y\rangle = \langle Ay,x\rangle$,
where
$\langle \cdot, \cdot \rangle$
is the symplectic structure. Alternatively, $A$ is self-adjoint in the Euclidean scalar product
$\langle J \cdot, \cdot \rangle$
and skew-commutes with $J$. The K\"ahler form on
$T_J \, Sp(2n,\R)/U(n)$
is given by
$\langle A,B\rangle = Tr \, JAB$.

\demo{3.3 Simplices on $Z$}
For two almost-complex structures
$J_1, J_2$,
tamed by $\omega$, we define a segment
$\Cal J(t) : t \mapsto (c(t,J_1(x),J_2(x))$
where
$c(t,J_1(x), J_2(x))$
is the geodesic segment in the Hermitian symmetric space of nonpositive curvature
$Sp(2n,\R)/U(n)$,
joining
$J_1(x)$
and
$J_2(x)$.
For a collection
$J_1, \ldots, J_n$
define a singular simplex
$K(J_1, \ldots, J_n)$
as in 1.3.
\enddemo
$   $\newline
{\it 3.4 Continuous cohomology classes in} $\text{Sympl}(M)$: {\it a definition}.
For any generator of the ring of
$Sp(2n,\R)$
-invariant form on
$Sp(2n,\R)/U(n)$
we define a continuous cohomology class in
$H_{cont}(\text{Sympl}(M),\R)$
by the explicit formula

$$\delta (f_1, \ldots, f_n) = \int_{K(J_0,{f_1}_* J_0 \ldots,{f_1 f_2 ...f_n}_* J_0 )} \varphi$$
where
$J_0$
is any fixed tamed almost-complex structure, and
$\varphi$
is a form of 3.2.

\demo{3.5 Non-trviality}
Let $M$ be a flat torus
$\R^{2n}/\Z^{2n}$
with a standard symplectic structure
$dx_1 \wedge dx_2 + \ldots + dx_{2n-1} \wedge dx_{2n}$.
As in 2.1, we have an
$Sp(2n,\Z)$
-invariant embedding
$Sp(2n,\R)/U(n) \hookrightarrow X$,
and the classes of 3.4 on
$\text{Sympl}(M)$
restrict to Borel classes on
$Sp(2n,\Z)$,
nontrivial for big $n$ [B].
\enddemo

\demo{\bf {3.6 Application to moduli spaces}}
Let $S$ be a closed Riemann surface of genus
$g \ge 2$,
and let
$\Cal M_g$
be a component of the representation variety
$Hom(\pi_1(S), SO(3))/SO(3)$
with Stiefel-Whitney class 1. This is known to be a smooth compact simply-connected symplectic manifold [Go2] of dimension
$6g-6$.
By a famous theorem of [NS],
$\Cal M_g$
is identified with the moduli space of stable holomoprhic vector bundles of rank 2 and odd determinant. The mapping class group
$Map_g$
acts symplectically on
$\Cal M_g$,
so we have an injective homomorphism
$Map_g \rightarrow Sympl(\Cal M_g)$.
Now we claim the following
\enddemo

\proclaim{Theorem (3.6)}
$H^2(Sympl(\Cal M_g),\R)$
is nontrivial. Morover, the homomorphism
$Map_g \rightarrow Sympl(\Cal M_g,\R)$
induces a nontrivial map in second real cohomology.
\endproclaim

\demo{Proof}
By the main theorem of [NS] there is a holomorphic embedding of the Teichm\"uller space
$T_g$
to the space of complex structures in
$\Cal M_g$,
tamed by Goldman's symplectic form. In particular, we have a
$Map_g$-invariant
holmorphic embedding
$T_g \overset{\alpha} \to \longrightarrow Z(\Cal M_g)$.
Let
$\Omega$
be the K\"ahler form of
$Z(\Cal M_g)$,
then
$\alpha^*(\Omega)$
is a 
$Map_g$
-equivariant K\"ahler form on
$T_g$.
We know there exist holomorphic maps
$Y \overset{\pi} \to \longrightarrow S$,
where $S$ is a closed Riemann surface, $Y$ is a compact complex surface and
$\pi$
is a smooth fibration by complex curves of genus $g$, such that the corresponding holomorphic map
$\tilde S \rightarrow T_g$
is nontrivial. We may form a flat holomorphic fibration
$\Cal F \rightarrow S$
with
$T_g$
as a fiber, associated to the homomorphism
$\pi_1(S) \rightarrow Map_g$,
coming from
$\pi$.
The Borel regulator of the flat fibration
$\Cal F \rightarrow S$,
corresponding to the form
$\alpha^*(\Omega)$
on $T_g$,
will coincide with the pullback of the class in
$H^2(Sympl(\Cal M_g),\R)$
under the composite map
$\pi_1(S) \rightarrow Map_g \rightarrow Sympl(\Cal M_g)$.
The variation of complex structure
$Y \overset \pi \to \longrightarrow S$
gives a holomorphic section of
$\Cal F \rightarrow S$
which is not horizontal. Therfore the pullbak of
$\alpha^*(\Omega)$
on $S$ using this section will have positive integral over $S$. By [Re1], section 3, this precisely means that the class we get in
$H^2(S,\R)$
is nontrivial. Therefore the map
$Map_g \rightarrow Sympl (\Cal M_g)$
induces a nontrivial map in
$H^2$.
\hfill Q.E.D.
\enddemo

\head{4. Bounded cohomology for area-preserving diffeomorphisms}\endhead

\demo{4.1}
Let 
$M^2$
be a compact oriented surface of any genus and let
$\nu$
be an area from on $M$. Then
$Diff_{\nu}M = \text{Sympl}(M)$.
The construction of 3.4 gives a class in
$H^2_{cont}(Diff_{\nu}M,\R)$.
\enddemo

\proclaim{Theorem (4.1)}
The cocyle
$\delta(h_1,h_2)$
of 3.4 is bounded. The class
$[\delta]$
lives therefore in the image of the natural map

$$H^2_b(Diff_{\nu}(M),\R) \rightarrow H^2(Diff^{\delta}_{\nu}(M),\R)$$
\endproclaim

\demo{Proof}
Fix a tame almost-complex structure
$J_0$.
Then
$\delta(h_1,h_2)$
is given by
$\int_M \, \text{area}_h(\qquad) \cdot \omega$,
where
$\text{area}_h(x,y,z)$
is the hyperbolic area in
$SL_2(\R)/SO(2) \approx \Cal H^2$
of the geodesic triangle, spanned by
$x,y,z$.
Therefore
$|\delta (h_1,h_2)| \le \pi \cdot \omega(M)$.
\enddemo

\demo{4.2 Non-amenability of two-generated subgroups of
$Diff_{\nu}(M)$}
We will apply theorem 4.1 to the following problem: given two area-preserving maps
$f,g : M \rightarrow M$,
when the group
$\phi(f,g) \in Diff_{\nu}(M)$
is ``big'' (say, free)? When
$Diff_{\nu}(M)$
is replaced by a finite-dimensional Lie group, this problem has been studied extensively, see e.g. [Re4], and references therein. In [Re4] we showed how the value of a (twisted) Euler class forces
$2 \kappa$
elements
$f_1, \ldots, f_{2\kappa}$
of
$SL_2(\R)$
to generate a free group. Here we will give a criterion for
$\phi(f,g)$
as above to be non-amenable. For that, denote
$F(f,g)$
a free group in two generators
$f,g$.
Consider the
$\ell^1$
-homology Banach space
$H^{\ell^1}_2(F,\R)$
[M]. An element of this space has a representive
$\sum^{\infty}_{j=1} a_j(h_j,k_j)$
with
$h_j,k_j \in F, \Sigma |a_j| < \infty$
and
$\sum \, a_j(h_jk_j -h_j-h_j)=0$
in
$\ell^1(F)$.
A bounded cocycle
$\ell$
induces a continuous functional

$$\sum \, a_j \, \ell(h_i,k_i) : H^{\ell_1}_2(F,\R) \rightarrow \R$$
which vanishes if
$[\ell] =0$
in
$H^2_b(F,\R)$.
\enddemo

\proclaim{Theorem (4.2)}
Let
$\sum \, a_j (h_j,k_j)$
be any
$\ell^1$
-cycle in
$H^{\ell_1}_2(F,\R)$.
If
$\sum \, a_j \, \delta(h_j,k_i) \ne 0$,
then the group
$\phi(f,g)$
is non-amenable. The set of pairs
$(f,g) \in Diff_{\nu}(M) \times Diff_{\nu}(M)$
satisfying this inequality, is open in
$C^1$
-topology.
\endproclaim

\demo{Proof}
Consider the following maps:

$$H^2_b(Diff_{\nu}(M),\R) \rightarrow H^2_b(\phi(f,g), \R) \rightarrow H^2_b(F(f,g), \R) \rightarrow (H^{\ell_1}_2(F(f,g),\R))^*$$
If
$\phi(f,g)$
is amenable, then
$H^2_b(\phi(f,g),\R) =0$
[Gr2], so the image of
$\delta$
in
$(H^{\ell_1}_2(F(f,g),\R))^*$
is zero and
$(\delta,\sum \, a_i(h_j,k_j))=0$,
a contradiction. The last statement of the theorem is checked directly from the definition of
$\delta$.
\enddemo

\demo{4.3 Constructing $\ell^1$-cycles}
The cardinality of
$\dim_{\R} H^{\ell_1}_2(F(f,g),\R)$
is
$2^{\aleph_0}$
by [M]. To apply the theorem 4.2 it is useful to have explicit formulas for 
$\ell^1$
-cycles. One way is described in [M].
\enddemo

\head{5. Lie algebra cohomology}\endhead

We will give the Lie algebraic analogues of the above constructed classes in
$Diff_{\nu}(M)$
and
$\text{Sympl}(M)$.
Observe that some odd-dimensional classes in the Lie algebra of
$Sympl(M)$
were constructed in [Re2] they induce, in general, nontrivial classes in cohomology of
$\text{Sympl}(M)$
as a topological space. The even-dimensional classes constructed here always induce trivial classes in
$H^*(\text{Sympl}^{top}(M),\R)$.

\demo{5.1 Formulas for $Diff_{\nu}(M)$}
Let
$X_1, \ldots, X_{2\kappa +1} \in \text{Lie} (Diff_{\nu}(M))$.
Fix a Riemannian metric $g$ with volume form $V$. Let

$$\psi(X_1,\ldots, X_{2\kappa+1}) = \int_M Alt \, Tr \prod^{2\kappa +1}_{j=1} (\bigtriangledown X_j+ (\bigtriangledown X_j)^*) \cdot \nu$$
\enddemo

\proclaim{Theorem (5.1)}
$\psi$
defines a cocycle for
$H^{2\kappa +1}(\text{Lie} (Diff_{\nu}(M))$.
\endproclaim

\demo{Proof}
Consider a
$Diff_{\nu}(M)$
-equivariant evaluation map
$Diff_{\nu}(M) \rightarrow M : f \mapsto (f^*)^{-1}(g)$.
Then the
$Diff_{\nu}(M)$
-invariant forms on $M$, constructed in 1.1 induce left-invariant closed forms on
$Diff_{\nu}(M)$,
whose restriction on
$T_e \, Diff_{\nu}(M)$
will be a Lie algebra cocycle. The derivative of the evaluation map
$\text{Lie}(Diff_{\nu}(M)) \rightarrow T_g \, M$
is given by
$X \mapsto \Cal L_Xg =g(\bigtriangledown X + (\bigtriangledown X)^* \cdot, \cdot)$.
Accounting the formula for Borel classes (see e.g. [Re3]), one arrives above-written formula for
$\psi$.
\enddemo

\demo{5.2 Formulas for $\text{Sympl}(M)$}
Let
$X_1,\ldots, X_{2\kappa} \in \text{Lie} (\text{Sympl}(M))$.
Fix a tame almost-complex structure $J$. Let

$$\varphi_{2\kappa} (X_1,\ldots, X_{2\kappa}) = \int_M Alt \, Tr \, J \cdot \prod^{2\kappa}_{j=1} \Cal L_{X_j} J \cdot \omega^n$$
\enddemo

\proclaim{Theorem (5.2)}
$\varphi$
defines a cocycle for
$H^{2\kappa} (Lie(Sympl (M))$.
\endproclaim

\demo{Proof}
Same as for 5.1.
\enddemo

\demo{5.3 Vanishing for $\varphi_2$ for flat torus}
\enddemo

\demo{Proposition (5.3)}
Let
$M = \R^{2n}/\Z^{2n}$
be a torus with standard symplectic structure. Then for any choice of a tame almost-complex structure, the cohomology class of
$\varphi_2$
in
$H^2(\text{Lie}(\text{Sympl}(M)),\R)$
is zero.
\enddemo

\demo{Proof}
The cohomology class of
$\varphi_2$
does not depend on the choice of $J$, since $X$ is connected. Choose $J$ to be the standard complex structure. We need to work on the formula for
$\varphi_2$.
Let $g$ be a metric, defined  by
$g(J \, \cdot, \cdot) = \omega$
(flat in our case). We then have
$\Cal L_X J= [\bigtriangledown X,J]$
since $g$ is K\"ahler and
$\bigtriangledown_X J =0$.
So

$$\varphi_2(X,Y) = \int_M Tr \, J([\bigtriangledown X,J] [\bigtriangledown Y,J]- [\bigtriangledown Y,J] [ \bigtriangledown X,J]]) \cdot \omega^n$$
Let $X$ be Hamiltonian, so that
$X = J \, grad \, f$.
Then
$\bigtriangledown X = J \, H_f$,
where
$H_f$
is the Hessian of $f$. If $Y$ is also Hamiltonian, say
$Y = J \, grad \, h$,
we have

$$\varphi_2(X,Y) =-\int_M Tr \, J[H_f,J] [H_h,J] \cdot \omega^n$$
Direct computation shows that the last expression is zero for flat torus. Now,
$\text{Lie}(\text{Sympl}(M))$
is a semidirect product of the ideal of Hamiltonian vector fields and an abelian subalgebra of constant vector fields, generated by (multivalued) linear Hamiltonians. Clearly,
$\varphi_2(X,Y)$
is zero for all choices for $X$ and $Y$.
\enddemo

\demo{5.4 Vanishing of $\varphi_2$ for a symplectic surface}
\enddemo

\proclaim{Proposition (5.4)}
Let
$(M,\omega)$
be a compact surface with a symplectic form. Then for any choice of a tame almost-complex structure, the cohomology class of
$\varphi_2$
in
$H^2(\text{Lie}(Ham(M),\R))$
is zero.
\endproclaim

\demo{Proof}
Let $g$ be as above. Again we have

$$\varphi_2(X,Y) = -\int_M \, Tr \, J[H_f,J][H_h,J] \cdot \omega$$
The proposition follows now from the following remarkable identity.
\enddemo

\proclaim{Theorem (5.4)}
On a compact Riemannian surface
$(M,g)$
the following identity holds:

\TagsOnRight
$$\int_M Tr \, J[H_f,J][H_h,J] \cdot d \, \text{area} = -\int K(g)\{f,h\} \cdot d \, \text{area}, \tag *$$
where
$K(g)$
is the curvature of $g$.
\endproclaim

\demo{Proof}
We were only able to prove this identity by a direct (very) long computation ([Re1]), which we will sketch here. Let
$g = e^{A(x,y)}(dx^2 + dy^2)$
in local conformal coordinates. Then
$\Gamma^x_{xx} = \frac 1 2 A_x, \Gamma^y_{yy} = \frac 1 2 A_y$,
$\Gamma^x_{xy} = \frac 1 2 A_y, \Gamma^y_{xy} = \frac 1 2 A_x$,
$\Gamma^y_{xx} = -\frac 1 2 A_y$,
$\Gamma^x_{yy}= -\frac 1 2 A_x$.
Next,
$H_f = \bigtriangledown (Grad \, f)$
and to the matrix of
$H_f$
is

$$\pmatrix \format \l &\qquad \l \\
e^{-a} f_{xx} + \frac 1 2 e^{-A} (A_yf_y - A_xf_x) & e^{-A} f_{xy} - \frac 1 2 e^{-A} (A_yf_x + A_xf_y) \\
\text{   } \\
e^{-A}f_{xy} - \frac 1 2 e^{-A} (A_yfx + A_xf_y) & e^{-A} f_{yy} + \frac 1 2 e^{-A} (A_x f_x - A_y f_y) \endpmatrix$$
and the same for $h$. Substituting to the left side of (*) one gets

$$\gather -2 \left[ \int (e^{-A} f_{xy} - \frac 1 2 e^{-A} (A_y f_x + A_xf_y)) \cdot (h_{xx} - h_{yy} + A_y h_y - A_x h_x)- \right. \\
\text{  } \\
\left. -\int (e^{-A} h_{xy} - \frac 1 2 (A_yh_x + A_xh_y)) (f_{xx} - f_{yy} + A_y f_y - A_x f_x)\right] dxdy \endgather$$
Twice integrating by parts, one finds this equal to

$$\gather \int e^{-A} [ -A_{xxy} f_x + A_y A_{xx} f_x - A_{yyy} f_x + \\
+ A_y A_{yy}f_x + A_{yyx} f_y -A_x A_{yy} f_y + A_{xxx} f_y - A_x A_{xx} f_y] dxdy \endgather$$
On the other hand, the right hand side is

$$\int_M \{ f_xh_y - f_yh_x\} \cdot (A_{xx} + A_{yy}) e^{-A} dx dy.$$
Again integrating by parts, one gets the same expression as above.
\hfill q.e.d.
\enddemo

\head{6. Chern-Simons-type class in $H^3(Diff_{\nu}(M^3),\R(\Z)$} \endhead

This section is best read in conjunction with [Re2]. In that paper, we constructed secondary classes in
$Hom (\pi_{2i-1} (B \, \text{Sympl}^{\delta} (M)^+, \R/A)$
where
$M^{2n}$
is a compact simply-connected symplectic manifold and $A$ is a group of periods of a biinvariant
$(2i-1)$
-form on
$\text{Sympl}(M)$,
whose restriction on the Lie algebra is
$f_1, \ldots, f_{2i-1} \rightarrow \text{Alt} \int_M \{f_1,f_2\} f_3 \ldots f_{2i-1} \cdot \omega^n$.
In particular, it implied the following results.

\proclaim{6.1 Theorem ([Re2]) (Chern-Simons class extends to $\text{Sympl}(S^2)$)}
There exists a rigid class in
$H^3(\text{Sympl}(S^2, can), \R/\Z)$
whose restriction on
$SO(3)$
is the standard Chern-Simons class.
\endproclaim

\proclaim{6.2 Theorem ([Re2]) (Chern-Simons class extends to $\text{Sympl}(\C P^2)$)}
There exists a rigid  class in
$H^3(\text{Sympl}(\C P^2, can), \R/\Z)$
whose restriction on 
$SU(3)$
is the standard Chern-Simons class.
\endproclaim

\proclaim{6.3 Theorem ([Re2])}
There exists a rigid  class in
$H^3(\text{Sympl}((S^2,a_1 \cdot can) \times S^2(a_2 \times can)), \R/\Z), a_1 \ne a_2$,
whose restriction on
$SO(3) \times SO(3)$
is the sum of standard Chern-Simons classes. 
\endproclaim

Let
$M^3$
be a rational homology sphere, say
$f \cdot H_1(M,\Z) =0, f \in \Z$.

\demo{6.4 The definition of the $ChS$ class}
Fix a point
$p \in M$
and consider the evaluation (at $p$) map

$$Diff_{\nu}(M) \rightarrow M.$$
The pull-back of
$\nu$
under this map is a closed left-invariant form
$\nu_p$
on
$Diff_{\nu}(M)$,
having integral periods. The general theory of [Re3] and [Re2] produces a regulator

$$\pi_3(B \, Diff^{\delta}_{\nu} (M)^+) \rightarrow \R/\Z \tag *$$
A different choice of a point
$p^{\prime} \in M$
will give another left-invariant form
$\nu_{p^{\prime}}$
such that
$\nu_p -\nu_{p^{\prime}} = d\mu$
for a left-invariant form
$\mu$.
It follows from [Re3] that the regulator (*) does not depend on $p$. In fact, one has a biinvariant 3-form $\omega$ on
$Diff_{\nu}(M)$,
whose values on the Lie algebra are given by
$\omega(X,Y,Z) = \int_M \nu(X(p), Y(p), Z(p))d \nu (p)$.
The form $\omega$ gives the same regulator as above.

To extend the regulator to
$H^3(Diff^{\delta}_{\nu}(M),\R/\Z)$,
we need to alter the scheme of [Re3] as follows. Since
$MSO_3(B \, Diff^{\delta}(M)) \approx H_3(B \, Diff^{\delta}(M),\Z)$
any class in
$H_3(B \, Diff^{\delta} (M), \Z)$
is represented by a map
$X \overset {\varphi} \to \longrightarrow B \, Diff (M)$,
or equivalently, by a representation
$\pi_1(X) \overset {\rho} \to \longrightarrow Diff_{\nu}(M)$.
Now, for $M$ a flat bundle
$M \rightarrow \Cal E \rightarrow X$,
associating to 
$\rho$.
The form $\omega$ extends to the closed form on
$\Cal E$
whose periods on fibers are 1. That gives an element
$\lambda$
in
$H^3(\Cal E, \R/\Z)$.
The spectral sequence of
$\Cal E$
with
$\R/\Z$
-coefficients looks like

$$\matrix \format \c & \qquad \c & \qquad \c & \qquad \c \\
\R/\Z & H^1(X,\R/\Z) & H^2(X,\R/\Z) & H^3(X,\R/\Z) \ldots \\
0 & 0 & 0 & 0 \\
H^0(X,\underline{W}) & H^1(X,\underline{W}) & H^2(X,\underline{W}) & H^3(X,\underline{W}) \ldots \\
\R/\Z & H^1(X,\R/\Z) & H^2(X,\R/\Z) & H^3(X,\R/\Z) \ldots
\endmatrix$$
where $W$ is the local system whose stalk at $p$ is
$H^1(M,\R/\Z) \approx \widehat{H_1(M,\Z)}$.
The element
$\lambda$
lies in the kernel of the wedge map
$H^3(\Cal E,\R/\Z) \rightarrow H^3(M,\R/\Z)$.
Now, the group
$H^2(X,\underline{W})$
has exponent a divisor of $f$, and the image of the transgression
$d^2 : H^1(X,\underline{W}) \rightarrow H^3(X,\R/\Z)$
has the same property. Therefore,
$f \cdot \lambda$
induces a well-defined class in
$H^3(X,\R/\Z \cdot \frac 1 f)$.
If $M$ is a
$\Z$-homology sphere, we get a class in
$H^3(X,\R/\Z)$.

If
$Y \rightarrow B \, Diff^{\delta}(M)$
is a map, bordant to
$\varphi$,
then the same argument as in [Re2] proves that the value of the corresponding class in
$H^3(Y,\R/\Z \cdot \frac 1 f)$
on $[Y]$ is the same as for $X$. So we constructed a well-defined map

$$H_3(Diff^{\delta}(M),\Z) \rightarrow \R/\Z \cdot \frac 1 f$$
\enddemo

\demo{6.5 Invariant scalar product on
$\text{Lie}(Diff_{\nu}(M)$,
the Cartan form and rigidity of ChS class}
Here we will prove that the ChS class

$$H_3(Diff^{\delta}_{\nu}(M),\Z) \rightarrow \R/\Z$$
of the previous section is rigid for
$M \approx S^3$.
For that purpose we need to work with principal flat bundles rather then with flat associated bundles. The clue is that the form $\omega$ constructed above on
$Diff_{\nu}(M)$
can be viewed as a Cartan form, associated with an invariant scalar product on
$\text{Lie}(Diff_{\nu}(M))$.
\enddemo

We are going to prove similar results for the group
$Diff_{\nu}(M^3)$
of volume-preserving diffeomorphisms of a compact oriented three-manifold. Throughout this section, $M$ is assumed to be a rational homology sphere, that is,
$H_1(M,\Z)$
is torsion.

Let
$X \in \text{Lie}(Diff_{\nu}(M))$
a vector field with
$\text{div} \, X =0$.
The form
$X \rfloor \nu$
is closed, whence exact:
$d\mu = X \rfloor \nu$.
Put
$\langle X,X\rangle = \int_M \mu \cdot (X \rfloor \nu)$.
An immediate computation shows that
$\langle X,X\rangle$
does not depend on the choice of
$\mu$.
Moreover
$X \mapsto \langle X,X\rangle$
is a quadratic form, invariant under the adjoint action of
$Diff_{\nu} (M)$.
By Arnold [A],
$\langle X,X\rangle$
is the asymptotic self-linking number of trajetories of $X$. We need the following elementary lemma (the proof of left to the reader)

\proclaim{Lemma (6.5)}
For any
$X,Y,Z \in \text{Lie} (Diff_{\nu}(M))$,

$$\Omega (X,Y,Z) = \omega(X,Y,Z)$$
that is, the forms 
$\Omega$
and $w$ coincide.
\endproclaim

Now, as in [Re2] we define a biinvariant form
$\Omega$
on
$Diff_{\nu}(M)$
by
$\Omega(X,Y,Z) = \langle [X,Y],Z\rangle$
on the Lie algebra.

\proclaim{Lemma (6.6)}
Let
$M = S^3/\Gamma$
where
$S^3$
is considered as a compact Lie group and the finite subgroup
$\Gamma$ acts from the right. Then the pullback of
$\Omega$
by the natural map
$S^3 \rightarrow Diff_{\nu}(M)$
is
$\frac 1 {|\Gamma |} \, \cdot$
(volume form of
$S^3$).
\endproclaim

\demo{Proof}
It is clearly enough to check this for
$\Gamma = \{1\}$.
Let
$v \in \text{Lie} (S^3)$
and $X$ is the corresponding right-invariant vector field. Let
$\mu$
be a right-invariant 1-form, defined by
$(v,\cdot)$
on Lie
$S^3)$.
Then
$d\mu = X \rfloor \nu$
and
$\mu \wedge (X \rfloor \nu) = \nu$.
\hfill q.e.d.
\enddemo

\proclaim{6.6 Theorem (Chern-Simons class in $Diff_{\nu} (S^3)$ )}
There exists a rigid class in
$H^3(Diff_{\nu}(S^3), \R/\Z)$
whose restriction on
$SO(4) \approx S^3 \times S^3/\Z_2$
coincides with the sum of standard Chern-Simons classes. Moreover, for
$M = S^3/\Gamma$
there exists a class in
$H^3(Diff_{\nu}(M),\R/\Z)$
whose restriction on
$S^3$
is 
$|\Gamma |$
times the standard Chern-Simons class.
\endproclaim

\demo{Proof}
By the general theory of regulators, developed in [Re3], section 3, and [Re2], the invariant form
$\Omega$
gives rise to a map

$$\pi_3 (B \, Diff^{\delta}_{\nu} (M)^+) \rightarrow \R/A$$
where $A$ is the group of periods of
$\Omega$
on the Hurewitz image of
$\pi_3(Diff_{\nu}(M))$
in
$H_3(Diff_{\nu}(M),\Z)$.
Moreover, if
$Diff_{\nu}(M)$
is homotopically equivalent to
$S^3$
or
$SO(4)$
this extends to a map

$$H_3(B \, Diff^{\delta}(M)) \rightarrow \R/A$$
By Hatcher [H] and Ivanov [I] this is exactly the case for
$M = S^3/\Gamma$.
Moreover, periods of
$\Omega$
are
$2\pi^2 \cdot \Z$
and
$2 \pi^2 \cdot \frac 1 {|\Gamma |} \Z$,
respectively. Since
$\Omega$
is a Cartan form, associated to an invariant polynomial in
$\text{Lie}(Diff_{\nu}(M))$,
it is rigid by Cheeger-Simons [Che-S].
\enddemo

\demo{6.6 Case of Seifert manifolds}
Let
$\Gamma$
be a uniform lattice in
$\widetilde{SL_2(\R)}$,
then
$M = \widetilde{SL_2(\R)}/\Gamma$
is a Seifert manifold. There is a cohomology class
$\beta \in H^3(\widetilde{SL_2(\R)}^{\delta},\R)$,
called the Seifert volume class [BGo], such that for any
$\Gamma \subset \widetilde{SL_2}(\R)$,
the restriction of 
$\beta$
on
$\Gamma$
is
$\text{vol} (\widetilde{SL_2(\R)}/\Gamma)$
times the fundamental class. Then the computation of 6.4 gives the class in
$H^3(Diff_{\nu}(M), \R)$,
whose restriction on
$\widetilde{SL_2(\R)}$
is
$\beta$,
subject to the condition that
$Diff_{\nu}(M)$
is contractible. It is not known to the author if this is true for all such $M$, comp. [FJ].
\enddemo

\head{7. Measurable transfer and higher asymptotic cycles}\endhead

We will first outline here an alternative approach in defining the classes of 1.2 in
$Diff_{\nu}(M)$.
For $M$ a locally symmetric space of nonpositive curvature, this approach also leads to new classes in
$H^*_{cont}(Diff_{\nu}(M),\R)$,
different from those of 1.2.

Let
$\frak G = Diff_{\nu}(M)$
and
$\frak G_0 \subset \frak G$
is a closed group, stabilizing a fixed point
$p \in M$.
Let
$\frak G^{\sim}$
be the connected component of
$\frak G$
and let
$\frak G^{\sim}_0 = \frak G \cap \frak G_0$.
Fix a measurable section
$s : M \rightarrow \frak G$
such that
$s(q)p = q$.
We will always assume that
$\overline{s(M)}$
is compact.

\demo{7.1 Ergodic cocycle in non-abelian cohomology [Gu]}
Define a map
$\psi : \frak G \times M \rightarrow \frak G_0$
by
$g \, s(q) = s(gq)\psi (g,q)$.
We will view it as a map
$\frak G \overset {\psi} \to \longrightarrow \Cal F(M,\frak G_0)$.
Here
$\Cal F(M,\frak G_0)$
is the group of measurable functions from $M$ to
$\frak G_0$
with compact closure of the image. 
$\frak G$
acts on
$\Cal F(M,\frak G_0)$
by the argument change and 
$\psi$
is a cocycle for the non-abelian cohomology
$H^1(\frak G,\Cal F(M,\frak G_0))$.
\enddemo

\demo{7.2 Measurable transfer [Gu]}
Now let
$f : \frak G_0 \times \ldots \frak G_0 \rightarrow \R$
be a locally bounded (say, continuous) cocycle. Define
$F : \frak G \times \ldots \times \frak G \rightarrow \R$
as
$F = \int_M f(\psi(g_1,m)$, 
$\psi(g_2,m) \ldots \psi(g_n,m))d \nu (m)$.
This defines a cohomology class in
$H^n(\frak G, \R)$,
independent of the choices of $s$ and $f$ [Gu].

Now, we have the tangential representation
$\frak G_0 \rightarrow SL(T_p(M))$.
Pulling back the usual Borel classes on
$\frak G_0$,
we construct cohomology classes in
$H^i(\frak G_0,\R)$
for
$i = 5,9, \ldots$.
The transfer will map these to classes in
$H^i(\frak G,\R)$,
which we have constructed in 1.2. We do not prove the comparison theorem here, however.
\enddemo

\demo{7.3 Supertransfer}
We will now define a map

$$H^{\kappa} (\pi_1(M),\R) \overset S \to \longrightarrow H^{\kappa}(Diff_{\nu}(M),\R)$$
in the following way. We know that
$\pi_0(\frak G^{\sim}_0) \approx \pi_1(M)/\pi_1(\frak G^{\sim})$.
This defines a homomorphism
$\frak G^{\sim}_0 \rightarrow \pi_0(\frak G^{\sim}_0) \rightarrow \pi_1(M)/\pi_1(\frak G^{\sim})$,
and a map
$H^{\kappa}(\pi_1(M)/\pi_1(\frak G^{\sim}),\R) \rightarrow H^{\kappa} (\frak G^{\sim}_0,\R)$.

In many interesting cases one knows that
$\pi_1(\frak G^{\sim}) = 1$.
If $M$ is a surface of genus
$g \ge 2$,
a result of Earle and Eells says that
$\frak G^{\sim}$
is contractible. For $M$ locally symmetric of rank
$\ge 2$
[FJ]. For any $M$ such that
$\pi_1(\frak G^{\sim}) =1$,
we get
$\pi_0(\frak G_0) \approx \pi_1(M)$
so that there is a map

$$H^k(\pi_1(M)) \rightarrow H^{\kappa} (\pi_0(\frak G^{\sim}_0)) \rightarrow H^{\kappa}(\frak G^{\sim}_0).$$
Now, composing with the measurable transfer
$H^{\kappa}(\frak G^{\sim}_0) \rightarrow H^{\kappa} (\frak G^{\sim})$
we arrive to a desired map

$$S : H^{\kappa} (\pi_1(M),\R) \rightarrow H^{\kappa} (\frak G^{\sim},\R)$$
\enddemo

\demo{7.4 Higher asymptotic cycles}
The dual to the above-constructed map $S$ is

$$S^{\vee} : H_{\kappa}(\frak G^{\sim},\R) \rightarrow H_{\kappa} (\pi_1(M),\R).$$
As we will see now, this is higher version of the classical asymptotic cycle character

$$\frak G^{\sim} \overset{\tau} \to \longrightarrow H_1(M,\R)$$
[Sch]. Indeed, for
$\kappa =1$
the map
$S^{\vee}$
will act as follows: let
$g \in \frak G^{\sim}$
be a volume-preserving map, isotopic to identity. Fix an isotopy
$g(t,x)$
such that
$g(0,\cdot) = \text{id}$
and
$g(1,\cdot) = g$.
For 
$x \in M, g(t,x)$
is a path from $x$ to
$g(x)$
and may be considered as a $1-$ current.
Now, the integral

$$\int_M [g(t,x)] d\nu (x)$$
is a closed current, defining an element in
$H_1(M,\R)$.
This will be
$S^{\vee} (g)$.

Now, the definition of the asymptotic cycle map [Sch] gives the following recepy: for an element
$z \in H^1(M,\Z)$
let
$f : M \rightarrow S^1$
be a representing map. The map
$f \circ g -f : M \rightarrow S^1$
is zero-homotopic, so it comes from the map
$F : M \rightarrow \R$.
Now,
$\int_MF(\text{mod} \, \Z)$
is the image of
$\tau(f)$
on $z$. If $f$ is isotopic to identity,
$\tau (f)$
lifts to
$H_1(M,\R)$.
It is easy to check that
$(df, \int_M[g(t,x)]d\nu) = (\tau (f),z)$,
which proves 
$S^{\vee} = \tau$
in dimension 1.
\enddemo


\enddocument